\documentclass[12pt]{article}

\renewcommand{\Re}{\textup{Re\,}}

\newcommand{\Ai}{\textup{Ai}}
\newtheorem{theorem}{Theorem}
\newtheorem{proposition}{Proposition}
\newtheorem{property}{Property}

\newtheorem{definition}{Definition}

\usepackage{amsmath}
\usepackage{amsfonts}
\usepackage{graphicx}

\title{Approximating the poles of a Painlev\'e V transcendent for the Fredholm determinant of the sine kernel}
\author{Walter Van Assche \\ KU Leuven, Belgium}
\date{\today}

\begin{document}
\maketitle

\begin{abstract}
The Fredholm determinant of the sine kernel has an expression as an integral containing a Painlev\'e transcendent $r$
satisfying Painlev\'e V, in the same way that the Fredholm determinant of the Airy kernel has an expression
involving the Hastings-McLeod solution of Painlev\'e II. In this paper we use Pad\'e approximation to approximate the poles of this Painlev\'e transcendent $r$ in the complex plane. We use Pad\'e approximation near $0$, Hermite-Pad\'e approximation near $0$ and two-point Pad\'e approximation near $0$ and infinity.
\end{abstract}

 \section{Introduction}
 The past few decades has seen an increased interest in non-linear special functions, in particular solutions of one of the Painlev\'e equations.
 It is known that these Painlev\'e equations have rational solutions for some particular values of the parameters appearing in the equations,
 and also solutions in terms of determinants containing special functions (Airy function, Bessel functions, hypergeometric and confluent hypergeometric functions, 
 parabolic cylinder functions). However most of the solutions will be new transcendental functions, which are called Painlev\'e transcendents.
 Some of these new transcendental functions appear in various applications, such as the Hastings-McLeod function, a special solution of Painlev\'e II which appears in the Tracy-Widom distribution, or the more general Ablowitz-Segur functions \cite{Bertola}.
Painlev\'e I has no parameters and all solutions are meromorphic functions with infinitely many poles. Boutroux \cite{Boutroux} 
showed that the poles are in sectors $S_k$ of the complex plane,
\[  S_k = \{ \frac{2k-1}{5} < \arg z < \frac{2k+1}{5} \}, \qquad 0 \leq k \leq 4, \]
and that for any two adjacent sectors there is a one-parameter solution of Painlev\'e I for which the poles (infinitely many) are in the remaining three sectors. These are called
tronqu\'ee solutions. For any four sectors there is also one solution with no poles in the four sectors and all poles (infinitely many) in the remaining sector. This is called a tritronqu\'ee solution. The Hastings-McLeod solution of Painlev\'e II has no poles in $\{ -\frac{\pi}{3} \leq \arg z \leq \frac{\pi}{3} \} \cup \{ \frac{2\pi}{3} \leq \arg z \leq \frac{4\pi}{3} \}$ and all the poles (infinitely many) are in the remaining sectors, see \cite{HXZ}. Tronqu\'ee solutions for Painlev\'e III and IV are investigated in \cite{Xia}.

In the present paper we make a numerical investigation of a special Painlev\'e transcendent for Painlev\'e V which appears when expressing the Fredholm determinant of the sine kernel. This special solution is free of poles in $\{ -\frac{\pi}{2} \leq \arg z \leq \frac{\pi}{2} \}$ and seems to have infinitely many poles in the other half plane $\Re z < 0$.
Our numerical experiments are based on various Pad\'e approximations, very much  in the spirit of Novokshenov \cite{Novok}.

In Section \ref{sec:Fred} we give a short explanation of Fredholm determinants, in particular for the sine kernel, the Airy kernel and the Bessel kernel.
These Fredholm determinants can be expressed in terms of a special Painlev\'e transcendent, e.g., the Airy kernel uses the Hastings-McLeod transcendent.
In Section \ref{sec:sine} we take a closer look at the sine kernel and we list what is known about the Painlev\'e transcendent in the expression of the Fredholm determinant.
Then in Section \ref{sec:Pade} we use Pad\'e approximation near $0$ and observe that Pad\'e approximation near infinity does not give any information about the poles of the
Painlev\'e transcendent. In Section \ref{sec:HermitePade} we use Hermite-Pad\'e approximation near zero,
and finally in Section \ref{sec:twopoint} we use two-point Pad\'e approximation near $0$ and infinity. In the final section we give a conclusion of our results.

\section{Fredholm determinants}   \label{sec:Fred}

Let $K(x,y)$ be a kernel on $\mathbb{R}^2$ and consider the operator $K_{[a,b]}$ on $L^2([a,b])$ for which
\[     K_{[a,b]} f (x) = \int_a^b K(x,y) f(y) \, dy.   \]

\begin{definition}
The Fredholm determinant of $K_{[a,b]}$ is given by
\[ \det ( I - K_{[a,b]}) = 1 + \sum_{n=1}^\infty \frac{(-1)^n}{n!} \int_a^b \cdots \int_a^b K \begin{pmatrix} x_1& \cdots & x_n \\ x_1 & \cdots & x_n \end{pmatrix}\,
   dx_1 \ldots dx_n,  \]
where
\[     K \begin{pmatrix} x_1& \cdots & x_n \\ x_1 & \cdots & x_n \end{pmatrix} = \det \Bigl( K(x_i,x_j) \Bigr)_{i,j=1}^n.   \]
\end{definition}
We assume that the infinite series converges; a sufficient condition is that $K_{[a,b]}$ is a trace class operator.

\subsection{The sine kernel}
The sine kernel is
\[     S(x,y) = \frac{\sin \pi(x-y)}{\pi(x-y)}. \]
Let $s<t$, then the operator $S_{[s,t]}$ maps $f  \in L^2([s,t])$ to $S_{[s,t]}f$ with
\[     S_{[s,t]}f(x) = \int_s^t S(x,y)f(y)\, dy.  \]
The Fredholm determinant is
\[  \det (I-S_{[s,t]}) = 1 + \sum_{n=1}^\infty \frac{(-1)^n}{n!}\int_s^t\cdots\int_s^t S\begin{pmatrix} x_1 & \cdots & x_n \\ x_1 & \cdots & x_n \end{pmatrix}\, dx_1\ldots dx_n, \]
where
\[    S\begin{pmatrix} x_1 & \cdots & x_n \\ x_1 & \cdots & x_n \end{pmatrix} = \det \Bigl( S(x_i,x_j) \Bigr)_{i,j=1}^n.  \]

For the point process with the sine kernel one has
\[    \textup{Prob}\bigl(\textup{no points in } [s,t]\bigr)  =  \det (I-S_{[s,t]}).   \] 

In random matrix theory one has for a random matrix $M_n$ from the Gaussian Unitary Ensemble (GUE): 
if $x^*$ is a point in the bulk of the spectrum $(-2,2)$, then
\begin{equation*}  
  \lim_{n \to \infty} \textup{Prob} \left( \textup{no eigenvalues of } M_n \textup{ in } \Bigl( x^*+\frac{s}{n\psi(x^*)},x^*+\frac{t}{n\psi(x^*)}\Bigr) \right) 
     = \det (I-S_{[s,t]}) ,  
\end{equation*}  
where $\psi$ is the density of the asymptotic eigenvalue distribution on $(-2,2)$. This result holds for a large class of unitary ensembles and is known as the gap probability in the bulk of the spectrum.

Jimbo, Miwa, M\^ori and Sato \cite{JMMS} found an explicit formula for this Fredholm determinant.
\begin{theorem}  \label{thm:Fsin}
The Fredholm determinant $\det(1-S_{[-\frac{x}{2},\frac{x}{2}]})$ is given by
\begin{equation}  \label{Fsin}
  \exp \left( -x - \int_0^x (x-t) r^2(t)\, dt \right)  
\end{equation}  
where $r(t)$ satisfies the differential equation
\begin{equation}  \label{rdiffeq}
 \bigl( (tr)''+\pi^2tr\bigr)^2 = 4r^2 \Bigl( (\pi tr)^2 + \bigl((tr)'\bigr)^2\Bigr),  
\end{equation} 
with initial conditions $r(0)=1$ and $r'(0)=1$.
\end{theorem}

The differential equation \eqref{rdiffeq}
is related to the sigma-form of Painlev\'e V
\[   ( t\sigma'')^2 - \bigl( 2(\sigma')^2 - t\sigma +\sigma\bigr)^2 + 4 (\sigma'+\kappa_1)(\sigma'+\kappa_2)(\sigma'+\kappa_3)(\sigma'+\kappa_4)=0.  \]

\subsection{The Airy kernel}
The Airy kernel is given by
\[   A(x,y) = \frac{\Ai(x)\Ai'(y)-\Ai'(x)\Ai(y)}{x-y},   \]
where $\Ai(x)$ is the Airy function, which satisfies $y''(x)=xy(x)$ and $\lim_{x \to \infty} \Ai(x) = 0$.
Consider the operator $A_{[s,\infty)}$ for which
\[   A_{[s,\infty)} f(x) = \int_s^\infty A(x,y) f(y)\, dy,   \]
then
\[ \det (I - A_{[s,\infty)}) = 1 + \sum_{n=1}^\infty \frac{(-1)^n}{n!} \int_s^\infty \!\! \cdots \int_s^\infty A \begin{pmatrix} x_1 & \cdots & x_n \\ x_1 & \cdots & x_n \end{pmatrix} \, dx_1\ldots dx_n, \]
where
\[   A \begin{pmatrix} x_1 & \cdots & x_n \\ x_1 & \cdots & x_n \end{pmatrix}  = \Bigl( A(x_i,x_j) \Bigr)_{i,j=1}^n.  \]
If $b$ is a \textbf{soft edge} of the spectrum where the eigenvalue density satisfies
\[   \psi(x) = \frac{c}{\pi} (b-x)^{1/2} \bigl( 1+ o(1) \bigr), \qquad x \to b, \]
then
\[   \lim_{n \to \infty}  \textup{Prob} \left( \textup{no eigenvalues of } M_n \textup{ in } \Bigl[ b+ \frac{s}{(cn)^{2/3}},\infty\Bigr) \right) = \det (1-A_{[s,\infty)}).  \]
For the largest eigenvalue $\lambda_n$ one has
\[  \lim_{n \to \infty}  \textup{Prob} \left(  c n^{2/3} (\lambda_n - b) < s \right) =  \det (1-A_{[s,\infty)}). \]

Tracy and Widom \cite{TW1} found an explicit formula, which is nowadays known as the Tracy-Widom distribution:

\begin{theorem}   \label{thm:FAiry}
\begin{equation}   \label{FAiry}
  \det(I-A_{[s,\infty)}) =    \exp \left( - \int_s^\infty (t-s)q^2(t)\, dt \right)  
\end{equation}  
where $q$ satisfies Painlev\'e II
\[  q''(t) = 2 q^3(t) + tq(t),   \vspace{-6pt} \]
with 
\[     \lim_{t \to \infty}  \frac{q(t)}{\Ai(t)} = 1.    \]   
\end{theorem}

This special solution of Painlev\'e II is known as the Hastings-McLeod solution. The distribution of the poles in the complex plane of the Hastings-McLeod solution was investigated by Huang et al. \cite{HXZ}. See also Bertola \cite{Bertola} and Bothner \cite{Bothner} for a larger family of solutions of Painlev\'e II.

\subsection{The Bessel kernel}
The Bessel kernel is ($\nu \geq 0$, $x,y \geq 0$)
\[   B(x,y) =     \frac12 \frac{J_\nu(\sqrt{x})\sqrt{y}J_\nu'(\sqrt{y})-\sqrt{x} J_\nu'(\sqrt{x})J_\nu(\sqrt{y})}{x-y},  \]
or equivalently
\[  B(x,y) = \frac12 \frac{J_\nu(\sqrt{x})\sqrt{y}J_{\nu-1}(\sqrt{y})-\sqrt{x} J_{\nu-1}(\sqrt{x})J_\nu(\sqrt{y})}{x-y}.  \]
This kernel appears as the limit of many kernels near a \textbf{hard edge}.

Let $B_t$ be the operator on $L^2([0,t])$ for which
\[   B_t f(x) = \int_0^t B(x,y) f(y)\, dy. \]
The Fredholm determinant is
\[ \det (1 - B_t) = 1 + \sum_{n=1}^\infty \frac{(-1)^n}{n!} \int_0^t \cdots \int_0^t  B \begin{pmatrix} x_1 & \cdots & x_n \\ x_1 & \cdots & x_n \end{pmatrix} \,
     dx_1 \ldots dx_n,   \]
where
\[     B \begin{pmatrix} x_1 & \cdots & x_n \\ x_1 & \cdots & x_n \end{pmatrix}  = \det \Bigl(  B(x_i,x_j) \Bigr)_{i,j=1}^n.   \]

If $b$ is a \textbf{hard edge} on the left of the spectrum where the eigenvalue density satisfies
\[   \psi(x) = \frac{c}{\pi} (x-b)^{-1/2} \bigl( 1+ o(1) \bigr), \qquad x \to b, \]
and $w(x)\sim (x-b)^\nu$ as $x \to b$ (from the right), then
\[   \lim_{n \to \infty}  \textup{Prob} \left( \textup{no eigenvalues of } M_n \textup{ in } \Bigl[b, b+ \frac{t}{4(cn)^{2}}\Bigr] \right) = \det (1-B_t).  \]

For the smallest eigenvalue $\lambda_1$ one has
\[   \lim_{n \to \infty} \textup{Prob} \left( 4(cn)^2 (\lambda_1-b) > t \right) = \det (1-B_t) .   \]

Tracy and Widom (again) \cite{TW2} 
gave an explicit formula:

\begin{theorem}   \label{thm:FBessel}
\begin{equation}    \label{FBessel}  
   \det (1-B_t) = \exp  \left( \frac{-1}{4}   \int_0^t (\log t/s ) p^2(s)\, ds \right), 
\end{equation}   
where $p$ satisfies the differential equation
\begin{equation*}               
   t(p^2-1)(tp')' = p(tp')^2 + \frac14 (t-\nu^2)p + \frac14 tp^3(p^2-2),       
\end{equation*}
with boundary condition
\[ p(t) \sim J_\nu(\sqrt{t}) \sim \frac{t^{\nu/2}}{2^\nu \Gamma(\nu+1)}, \qquad t \to 0.  \] 
\end{theorem}

The differential equation
\[    t(p^2-1)(tp')' = p(tp')^2 + \frac14 (t-\nu^2)p + \frac14 tp^3(p^2-2),   \]
is basically a version of Painlev\'e V.
Let
\[       p(t) = \frac{1+y(x)}{1-y(x)}, \qquad t=x^2, \]
then $y$ satisfies Painlev\'e V
\[   \frac{d^2y}{dx^2} = \left( \frac{1}{2y} + \frac{1}{y-1} \right) \left( \frac{dy}{dx}\right)^2 - \frac{1}{x} \frac{dy}{dx}
     + \frac{(y-1)^2}{x^2} \left( \alpha y + \frac{\beta}{y} \right) + \frac{\gamma y}{x} + \frac{\delta y(y+1)}{y-1}, \]
with $\alpha=\beta=\nu^2/8$, $\gamma = 0$ and $\delta=-2$. Note that $\gamma=0$ implies that this can be transformed to
Painlev\'e III, see \cite[\S 41]{GLS}.

\section{Back to the sine kernel}  \label{sec:sine}
\begin{figure}[ht]
\centering
\includegraphics[width=4in]{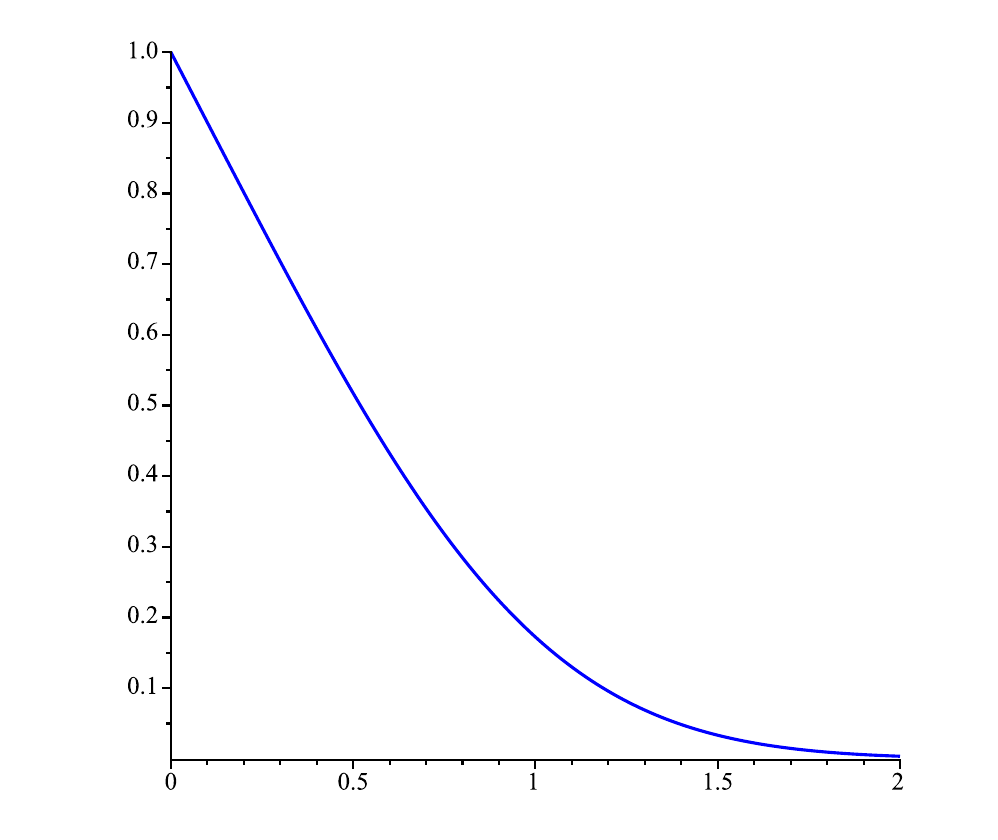}
\caption{The Fredholm determinant $F$ for the sine kernel}
\label{fig:Fsin}
\end{figure}

Recall
\[     S(x,y) = \frac{\sin \pi(x-y)}{\pi(x-y)}, \]
%
and Theorem \ref{thm:Fsin} giving an expression for the Fredholm determinant $F(x) = \det (1-S_{[-\frac{x}2,\frac{x}2]})$,
see Figure \ref{fig:Fsin}.

We will investigate the function $r$ in more detail.
Let $R_{[s,t]}$ be the operator for which
\[    I + R_{[s,t]} = (I - S_{[s,t]})^{-1},   \]
and $R(x,y;s,t)$ the corresponding kernel
\[        R_{[s,t]} f(x) = \int_s^t R(x,y;s,t) f(y)\, dy.   \]
For the analysis and proof of Theorem \ref{thm:Fsin} we refer to the original paper \cite{JMMS} or to \cite[\S 3.6]{AGZ}.
Here we only give a sketch and the ingredients we need for the rest of this paper.

\begin{proposition}
The function $r$ in \eqref{Fsin} is given by
\[   r(t) = R \left( -\frac{t}{2}, \frac{t}{2};-\frac{t}{2}, \frac{t}{2} \right).  \]
\end{proposition}

To find the differential equation for $r$, one needs to find the partial derivatives of $R(x,y;s,t)$:
\[     \frac{\partial R}{\partial x}, \quad  \frac{\partial R}{\partial y},   \quad  \frac{\partial R}{\partial s}, \quad  \frac{\partial R}{\partial t}.  \]
The latter two are in terms of two auxiliary functions $p$ and $q$
\[  p(t) = \cos \frac{\pi t}{2} + \int_{-t/2}^{t/2} \cos \pi y\ R\bigl(-\frac{t}{2},y;-\frac{t}{2},\frac{t}{2}\bigr)\, dy, \]
\[  q(t) = \sin \frac{\pi t}{2} + \int_{-t/2}^{t/2} \sin \pi y\ R\bigl(-\frac{t}{2},y;-\frac{t}{2},\frac{t}{2}\bigr)\, dy. \]

\subsection{The functions $p$ and $q$}
\begin{property}
\[    r(t) = \frac{2p(t)q(t)}{\pi t}, \]
and $(p,q)$ satisfy the system of first order differential equations
\begin{eqnarray}
    p' &=& -\frac{\pi}{2} q + \frac{2p^2q}{\pi t},   \label{p'} \\
    q' &=& \frac{\pi}{2} p - \frac{2pq^2}{\pi t}, \label{q'}
\end{eqnarray}
with $p(0)=1$, $q(0)=0$ and $q'(0)=\pi/2$.
\end{property}

\begin{figure}[ht]
\centering
\includegraphics[width=3in]{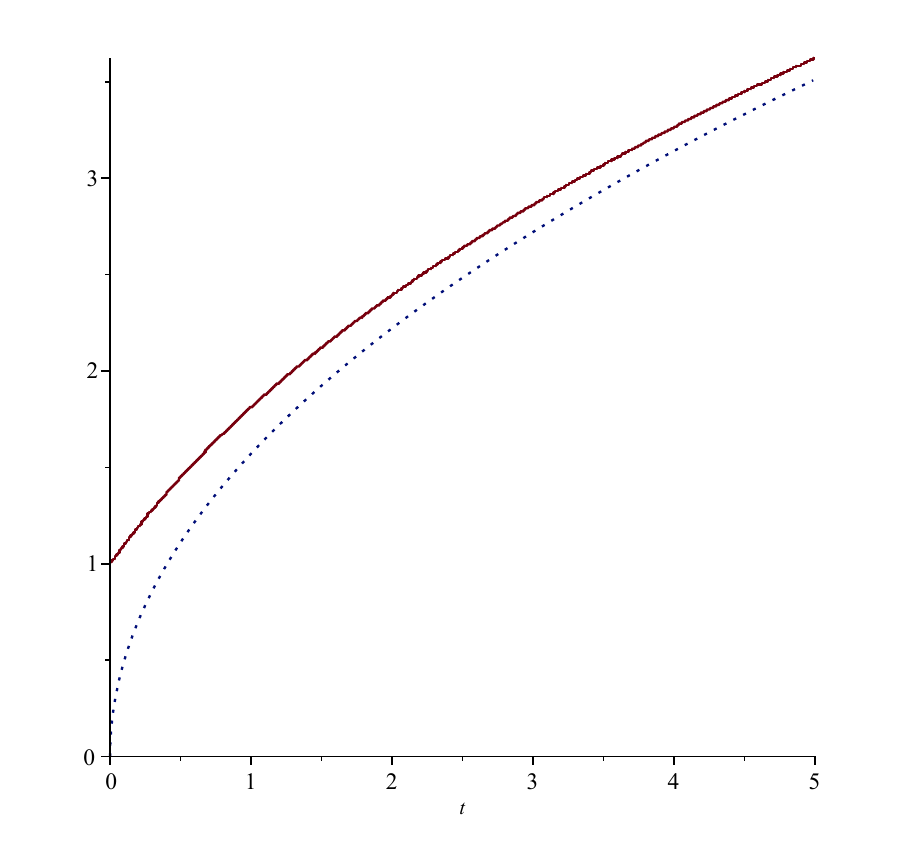}
\includegraphics[width=3in]{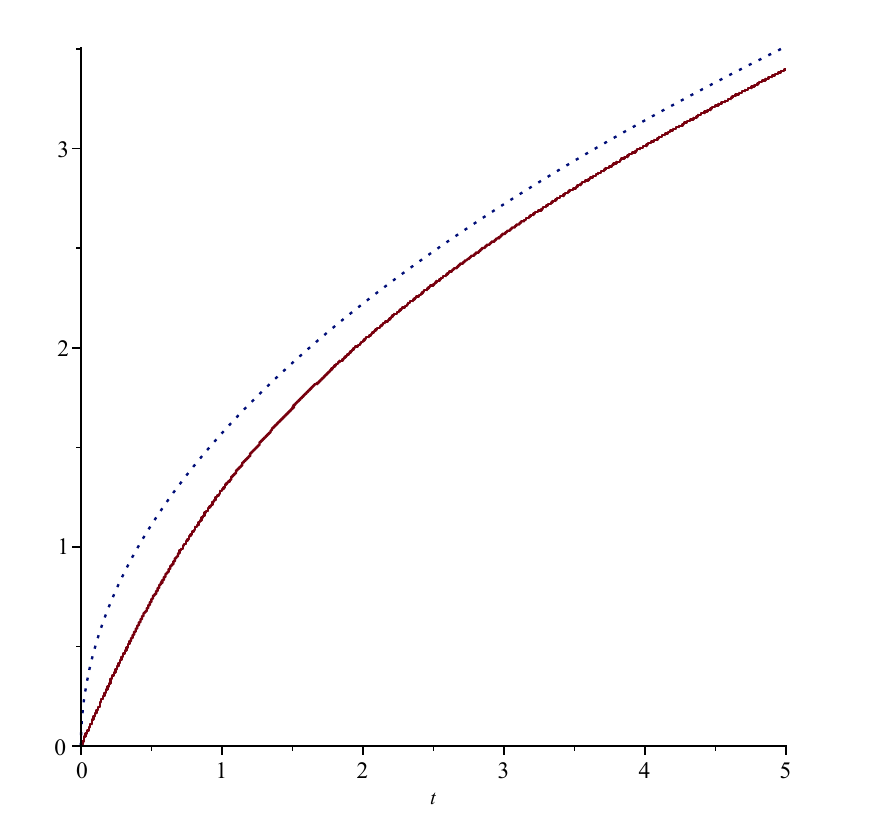}
\caption{the functions $p$ (left) and $q$ (right).}
\label{fig:pq}
\end{figure}

Both $p$ and $q$ satisfy a version of Painlev\'e V, 
\begin{align*}
   p'' &= (p')^2 \frac{p}{p^2-\frac{\pi^2}{4} t} - \frac{p'}{t} \frac{p^2}{p^2-\frac{\pi^2}{4}t} + \frac{p}{t} \bigl( p^2 - \frac{\pi^2}{4} t\bigr), \\
   q'' &= (q')^2 \frac{q}{q^2-\frac{\pi^2}{4} t} - \frac{q'}{t} \frac{q^2}{q^2-\frac{\pi^2}{4}t} + \frac{q}{t} \bigl( q^2 - \frac{\pi^2}{4} t\bigr).
\end{align*} 
The quantity
\[   \hat{\sigma} = -t(p^2+q^2) + \frac{4}{\pi^2} p^2q^2  \]
satisfies
\[ (t\hat{\sigma}'')^2 + 4\pi^2 (t\hat{\sigma}' -\hat{\sigma})\bigl(t\hat{\sigma}' -\hat{\sigma} +(\hat{\sigma}'/\pi)^2\bigr), \]
which relates to the $\sigma$-form of Painlev\'e V
\[   ( t\sigma'')^2 - \bigl( 2(\sigma')^2 - t\sigma' +\sigma\bigr)^2 + 4 (\sigma'+\kappa_1)(\sigma'+\kappa_2)(\sigma'+\kappa_3)(\sigma'+\kappa_4)=0.  \]
Indeed,  $\hat{\sigma}(t) = \sigma(2\pi it)$, for $\sigma$ in the standard form, with $\kappa_1=\kappa_2=\kappa_3=\kappa_4=0$.

Another expression for the Fredholm determinant is in terms of $\hat{\sigma}$:
\begin{equation}  \label{Fsigma}
   F(x) = \exp \left( \int_0^x \frac{\hat{\sigma}(t)}{t}\, dt \right). 
\end{equation}   
This follows from
\[   - \frac{\hat{\sigma}(t)}{t} = p^2+q^2-tr^2, \qquad   \frac{d}{dt} \frac{\hat{\sigma}(t)}{t} = -r^2.   \]

\begin{property}
The expression
\[   H(p,q) = \frac{\pi\hat{\sigma}}{4t} \]
is the Hamiltonian of this integrable system
\begin{eqnarray*}
  \frac{\partial H}{\partial p} = -\frac{\pi}{2} p + \frac{2pq^2}{\pi t} &=& - q'(t), \\
   \frac{\partial H}{\partial q} = \frac{\pi}{2} q + \frac{2p^2q}{\pi t} &=& p'(t).  
\end{eqnarray*}
\end{property}

\subsection{The function $r$}
Recall the differential equation \eqref{rdiffeq}
\[ \bigl( (tr)''+\pi^2tr\bigr)^2 = 4r^2 \Bigl( (\pi tr)^2 + \bigl((tr)'\bigr)^2\Bigr),  \]
with initial conditions $r(0)=1$ and $r'(0)=1$.

\begin{figure}[ht]
\centering
\includegraphics[width=4in]{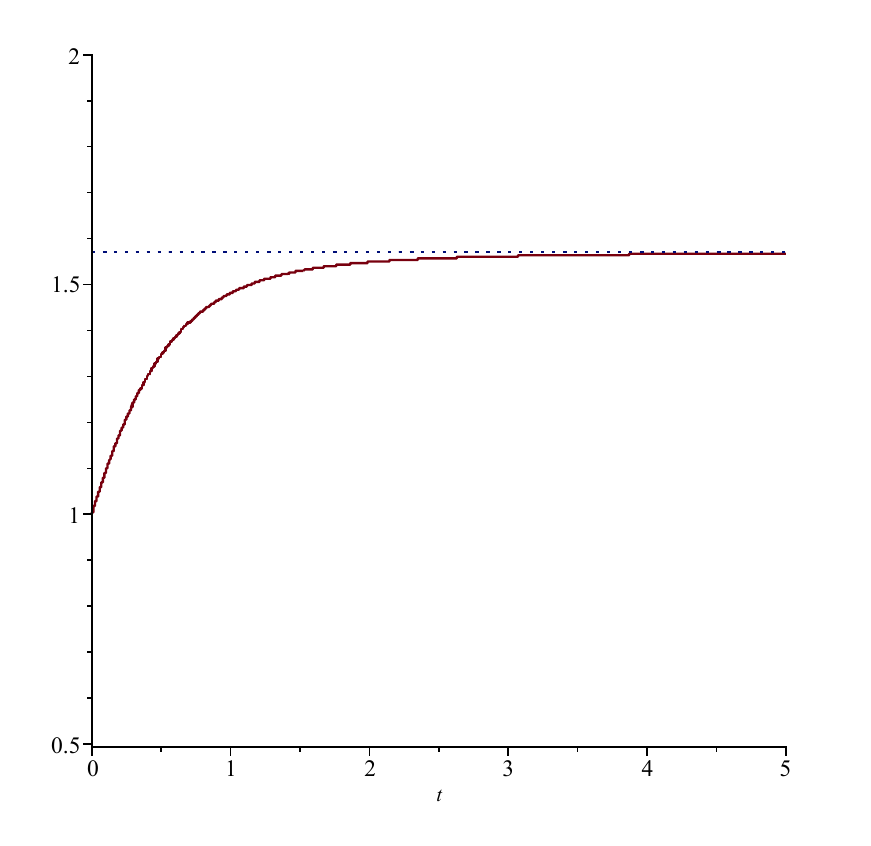}
\caption{The function $r$}
\label{fig:r}
\end{figure}

It is easier to use $r(t) = 2p(t)q(t)/\pi t$, 
where $(p,q)$ satisfy the system of first order differential equations
\begin{equation}  \label{pq}
     \begin{cases}   \   p' =  q \displaystyle \left( -\frac{\pi}{2}  + \frac{2p^2}{\pi t} \right),    \\[10pt]
                                   q' = p \displaystyle \left( \frac{\pi}{2}  - \frac{2q^2}{\pi t} \right),  \end{cases}  
\end{equation}
with initial conditions $p(0)=1$, $q(0)=0$, $q'(0)=\pi/2$.

\begin{property}
 The solution of the differential equation \eqref{rdiffeq} with $r(0)=1$ and $r'(0)=1$ has no poles on the
 positive real line and has the asymptotic behaviour
 \[    \lim_{t \to \infty} r(t) = \frac{\pi}{2}.   \]
 This corresponds to
 \[    p(t) \sim \frac{\pi}{2} \sqrt{t}, \quad q(t) \sim   \frac{\pi}{2} \sqrt{t}, \qquad t \to \infty.   \] 
\end{property}

See \cite{AGZ} or \cite{JMMS}.
This special solution is very sensitive to the choice of initial values.
\begin{itemize}
\item If $q(0) \neq 0$ then $p$ has a pole at $0$.
\item If $p(0) = 1+\epsilon$ with $\epsilon >0$, then $p$ has a pole at $t_0 >0$ where $q(t_0)=\frac{\pi}{2}\sqrt{t_0}$.
See Figure \ref{fig:p0=1.0001}.
\begin{figure}[ht]
\centering
\includegraphics[width=3in]{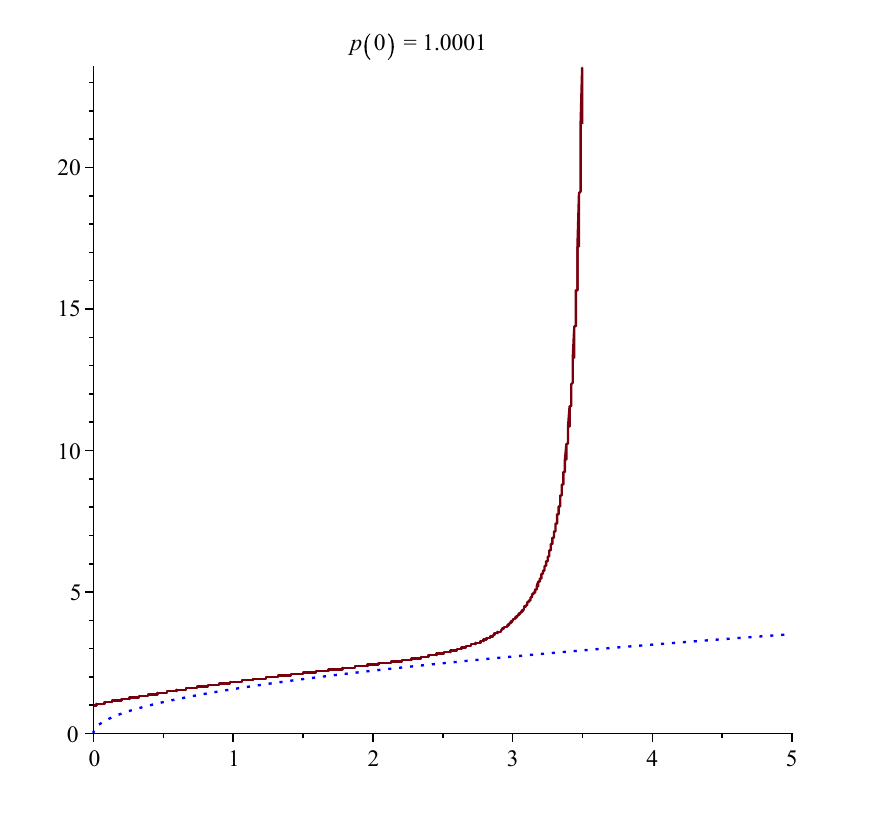}
\includegraphics[width=3in]{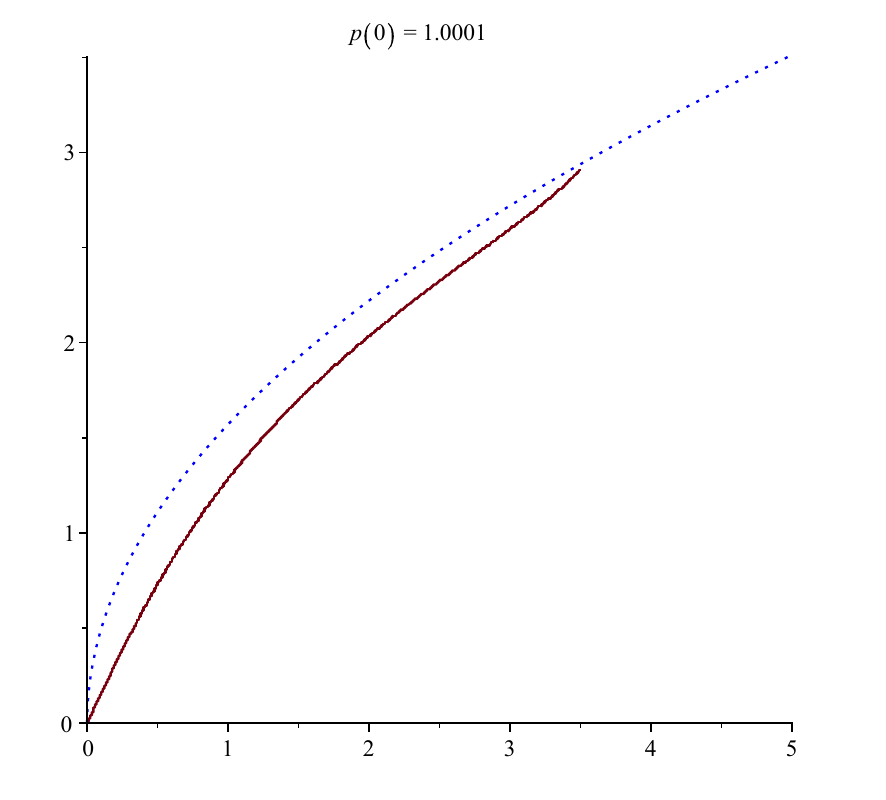}
\caption{The function $p$ (left) and $q$ (right) for $p(0)=1.0001$}
\label{fig:p0=1.0001}
\end{figure}
\item If $p(0) = 1-\epsilon$ with $\epsilon >0$ then $p$  crosses $\frac{\pi}{2}\sqrt{t}$ and oscillates around the horizontal axis.
See Figure \ref{fig:p0=0.999}.
\begin{figure}[ht]
\centering
\includegraphics[width=3in]{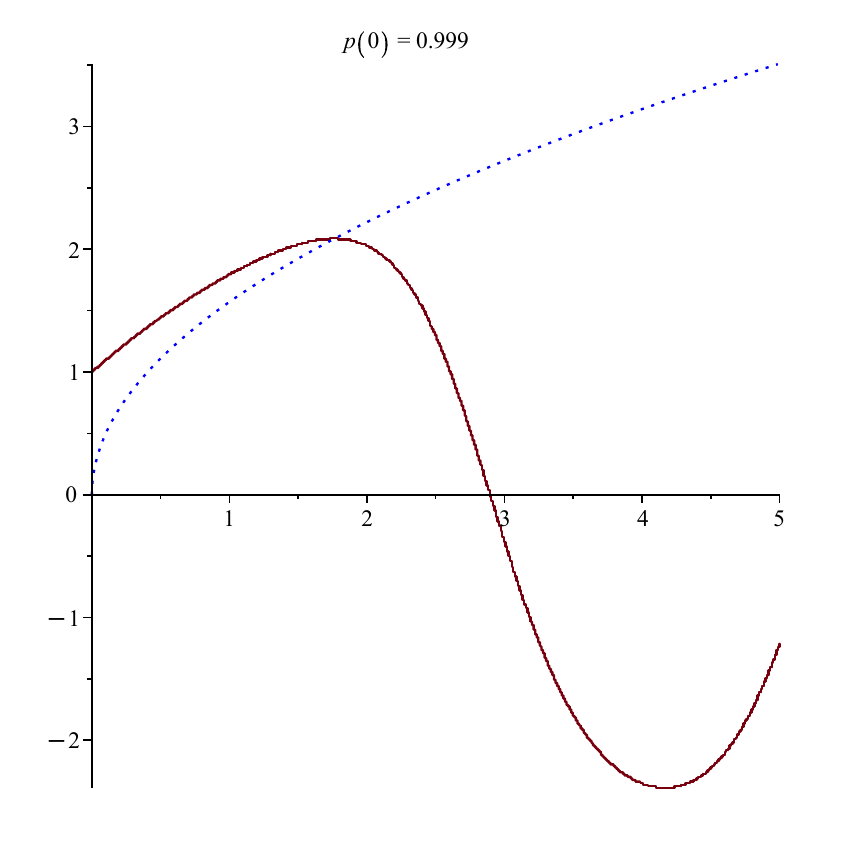}
\includegraphics[width=3in]{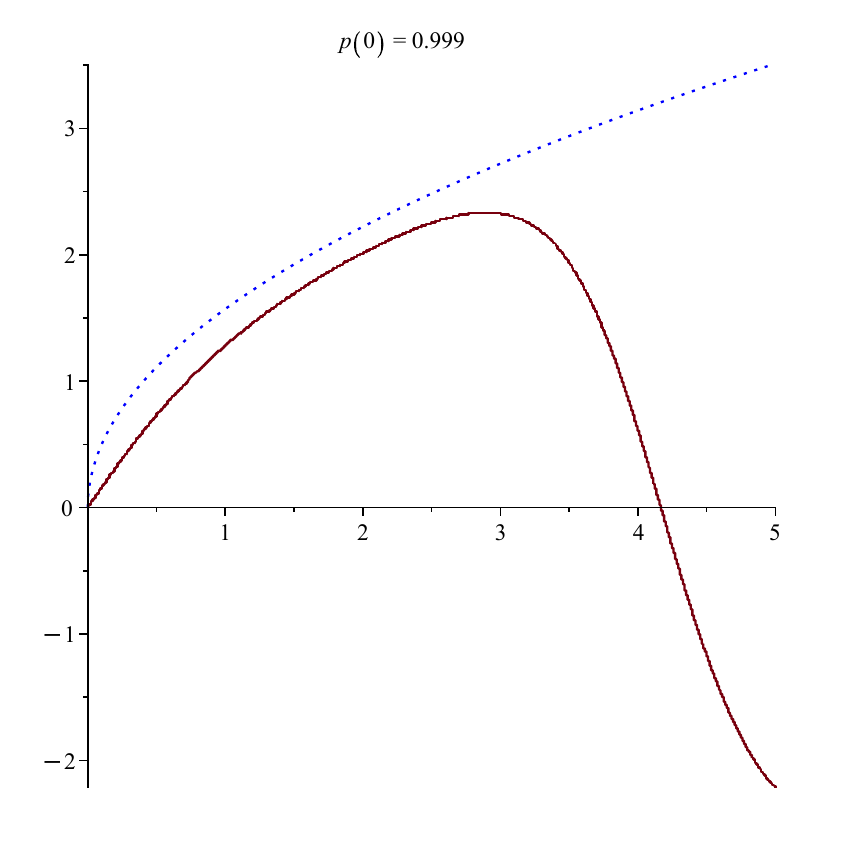}
\caption{The functions $p$ (left) and $q$ (right) for $p(0)=0.999$}
\label{fig:p0=0.999}
\end{figure}
\end{itemize}

Note that $p_0(t) = \frac{\pi}{2} \sqrt{t}$ and $q_0(t) =\frac{\pi}{2} \sqrt{t}$ are such that the right hand side
of \eqref{pq} vanishes (but not the left hand side), so that the function $t \mapsto \frac{\pi}{2}\sqrt{t}$  is expected to play a role. This function
is indicated with a dotted line in Figures \ref{fig:pq}, \ref{fig:p0=1.0001} and \ref{fig:p0=0.999}.

\subsection{Series near $t=0$}
If we write a series near $0$ as
\[   p(t) = \sum_{k=0}^\infty p_k t^k, \quad  q(t) = \sum_{k=0}^\infty q_k t^k, \]
and then insert this into \eqref{pq}, with initial conditions $p_0=1$ and $q_0=0$, then we get a non-linear  recurrence
\begin{align}
  (k+1)q_{k+1} &= \frac{\pi}{2} p_k - \frac{2}{\pi} \sum_{j=0}^{k+1} p_{k-j+1} \sum_{\ell=0}^j q_\ell q_{j-\ell}, 
       \label{recq0} \\
  (k+1)p_{k+1} &= -\frac{\pi}{2} q_k + \frac{2}{\pi} \sum_{j=0}^{k+1} q_{k-j+1} \sum_{\ell=0}^j p_\ell p_{j-\ell},
        \label{recp0}
 \end{align}
 The first few terms give
 \begin{multline*}
    p(t) = 1 + t + \left( 1-\frac{\pi^2}{8} \right) t^2 + \left( 1 - \frac{7\pi^2}{72} \right) t^3 + \left( 1-\frac{\pi^2}{8} + \frac{\pi^4}{384} \right) t^4 \\
     +\ \left( 1 - \frac{11\pi^2}{72} + \frac{143\pi^4}{28800} \right) t^5 + \left( 1 - \frac{13\pi^2}{72} +\frac{2171\pi^4}{259200} - \frac{\pi^6}{46080} \right) t^6 +\cdots 
\end{multline*}
\[  q(t) = \frac{\pi}{2} t - \frac{\pi^3}{48} t^3 + \frac{\pi^3}{72} t^4 + \frac{\pi^5}{3840} t^5 - \frac{13\pi^5}{14400} t^6
     + \left( \frac{\pi^5}{2592} - \frac{\pi^7}{645120} \right) t^7 + \cdots   \]     

\subsection{Series near infinity}   \label{sec:2.4}
In a similar way we can write a series near infinity
\[   p(t) = \frac{\pi}{2} \sqrt{t} \sum_{k=0}^\infty \hat{p}_k t^{-k}, \quad 
     q(t) = \frac{\pi}{2} \sqrt{t} \sum_{k=0}^\infty \hat{q}_k t^{-k}, \]
and if we insert this into \eqref{pq} with initial conditions $\hat{p}_0=1$ and $\hat{q}_0 = 1$, then we get the recursion
\begin{align*}
  2\pi \hat{p}_{k+1} &= (1-2k) \hat{p}_k 
  - \pi \sum_{j=0}^{k-1} \hat{q}_{k-j} \sum_{\ell=0}^{j+1} \hat{p}_{\ell}\hat{p}_{j+1-\ell}
  - \pi \sum_{\ell=1}^k \hat{p}_{\ell} \hat{p}_{k+1-\ell}, \\
   2\pi \hat{q}_{k+1} &= (2k-1) \hat{q}_k 
  - \pi \sum_{j=0}^{k-1} \hat{p}_{k-j} \sum_{\ell=0}^{j+1} \hat{q}_{\ell}\hat{q}_{j+1-\ell}
  - \pi \sum_{\ell=1}^k \hat{q}_{\ell} \hat{q}_{k+1-\ell}.
\end{align*}
The first few terms give  
\[   p(t) = \frac{\pi}{2} \sqrt{t} \left( 1+ \frac{1}{2\pi t} - \frac{1}{8\pi^2t^2} + \frac{5}{16\pi^3 t^3} 
  - \frac{85}{128 \pi^4 t^4} + \frac{719}{256 \pi^5 t^5} - \frac{11813}{1024\pi^6t^6} + \cdots \right) \]
\[  q(t) = \frac{\pi}{2} \sqrt{t} \left( 1 - \frac{1}{2\pi t} - \frac{1}{8\pi^2t^2} - \frac{5}{16\pi^3t^3} 
  - \frac{85}{128\pi^4t^4} - \frac{719}{256\pi^5t^5} - \frac{11813}{1024\pi^6 t^6} - \cdots \right).  \]
  Observe that $\hat{p}_k = (-1)^k \hat{q}_k$.
   
\subsection{Symmetric form of Painlev\'e V}
Put
\[    f_1(t) = - \frac{1}{\sqrt{\pi}} p(t/\pi) + \frac{\sqrt{t}}2, \quad  f_2(t) = \frac{1}{\sqrt{\pi}} q(t/\pi) + \frac{\sqrt{t}}2, \]
\[    f_3(t) =  \frac{1}{\sqrt{\pi}} p(t/\pi) + \frac{\sqrt{t}}2, \quad  f_4(t) = -\frac{1}{\sqrt{\pi}} q(t/\pi) + \frac{\sqrt{t}}2, \]
and $\mu_1=\mu_2=\mu_3=\mu_4=\frac14$, then
\begin{eqnarray*}
      tf_1' &=& f_1f_3(f_2-f_4)+(\frac12-\mu_3)f_1+\mu_1f_3 , \\
      tf_2' &=& f_2f_4(f_3-f_1)+(\frac12-\mu_4)f_2+\mu_2f_4 , \\
      tf_3' &=& f_3f_1(f_4-f_2)+(\frac12-\mu_1)f_3+\mu_3f_1 , \\
      tf_4' &=& f_4f_2(f_1-f_3)+(\frac12-\mu_2)f_4+\mu_4f_2 ,
\end{eqnarray*}
which is the symmetric form of Painlev\'e V \cite{NY}.
Then 
\[  w(t)=1-\frac{\sqrt{t}}{f_1(t)} = \frac{2p(t/\pi)+\sqrt{t\pi}}{2p(t/\pi)-\sqrt{t\pi}}  \]
satisfies  Painlev\'e V
\begin{equation*} 
w'' = \left( \frac{1}{2w} + \frac{1}{w-1} \right) (w')^2 - \frac{1}{t} w' + \frac{(w-1)^2}{t^2}\left( \alpha w + \frac{\beta}{w} \right) 
+ \gamma \frac{w}{t} + \delta \frac{w(w+1)}{w-1}
\end{equation*}
with $\alpha=-\beta=1/32$, $\gamma=0$ and $\delta=-1/2$.

\section{Pad\'e approximation}    \label{sec:Pade}
In this section we investigate the poles of the function $r$ in the complex plane by using Pad\'e approximants. See Novokshenov \cite{Novok} for this approach.

\subsection{Pad\'e approximation near $0$}   \label{sec:PA0}
The $[m/n]$ Pad\'e approximant of $r$ near $0$ is the rational function $Q_m/P_n$, with $\deg P_n \leq n$ and $\deg Q_m \leq m$, for which
\[    P_n(x) r(x) - Q_m(x) = \mathcal{O}(x^{m+n+1}), \qquad x \to 0.  \]
We can use the series expansions near $0$ to find the $[m/n]$ Pad\'e approximants of the function $r$ near $0$.
This allows to get an idea of the poles of $r$ near $0$ by using a theorem of de Montessus de Ballore, see
\cite[\S 6.2]{BG-M}

\begin{theorem}[de Montessus de Ballore]
Let $f$ be a function which is meromorphic in the disk $|z| \leq R$, with precisely $M$ simple poles at distinct points
$z_1,z_2,\ldots, z_M$, with $0 < |z_1| \leq |z_2| \leq \cdots \leq |z_M| < R$. Then
\[   \lim_{L \to \infty} [L/M](z) = f(z), \]
uniformly on any compact subset of $\mathcal{D}_M = \{ |z| \leq R, z \neq z_i, i=1,2,\ldots,M\}$.
\end{theorem}

\begin{figure}[ht]
\centering
\includegraphics[width=4in]{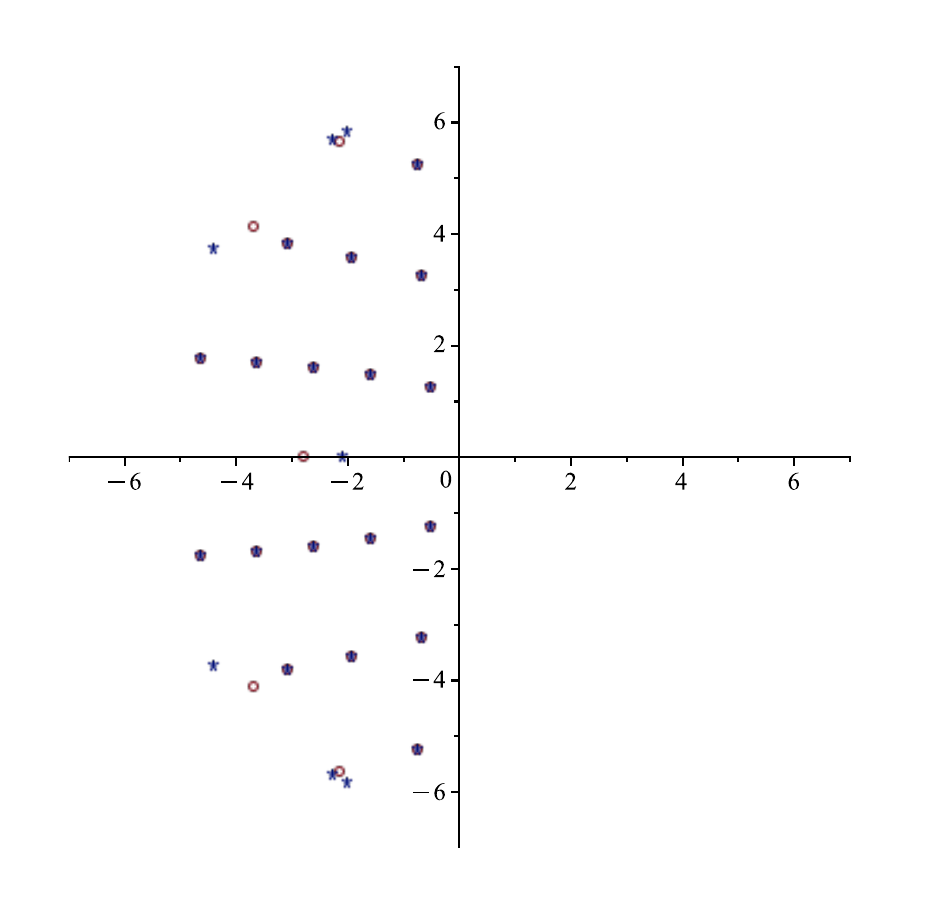}
\caption{The poles of the Pad\'e approximants $[50,25]$ (circles) and $[70,25]$ (asterisks).}
\label{fig:pade25}
\end{figure}

There is a more quantitative version of this theorem (due to Gonchar) that says how fast the poles of $f$ are approximated by the zeros of the denominator.
Suppose the $[L/M]$ Pad\'e approximant is $Q_L(x)/P_M(x)$ and let $\pi_M(x)$ be the polynomial that has the poles $z_1,\ldots,z_M$ as its zeros, then
Gonchar's theorem \cite[p. 135]{CYL2} says that
\[   \limsup_{L \to \infty} \| P_M-\pi_M \|^{1/L} = \frac{\max_{1 \leq k \leq M} |z_k|}{R} < 1,  \]
where $\|\cdot\|$ is the coefficient norm in the space of polynomials.

We computed the coefficients $p_k$ and $q_k$ numerically using the recurrence \eqref{recq0}--\eqref{recp0}, and set
\[     P(t) = \sum_{k=0}^{101} p_k t^k \approx p(t), \quad   Q(t) = \sum_{k=0}^{100} q_{k+1} t^k \approx q(t)/t, \] and then
used the command \texttt{pade} in Maple (package \texttt{numapprox} and \texttt{Digits:=100}) to compute
the Pad\'e approximants $[50/25]$ and $[70/25]$ for $R(t) = 2P(t)Q(t)/\pi$ near $t=0$. The poles of $[50/25]$ are plotted in Figure \ref{fig:pade25} as circles ($\circ$) and those of $[70/25]$ as asterisks ($\ast$). 
Observe that all the poles are in the half plane $\Re z <0$ and that the first 18 poles of both Pad\'e approximants are very close together. With de Montessus de Ballore
in mind, these 18 poles are expected to approximate 18 poles of the function $r$ closest to the origin.

\begin{figure}[ht]
\centering
\includegraphics[width=4in]{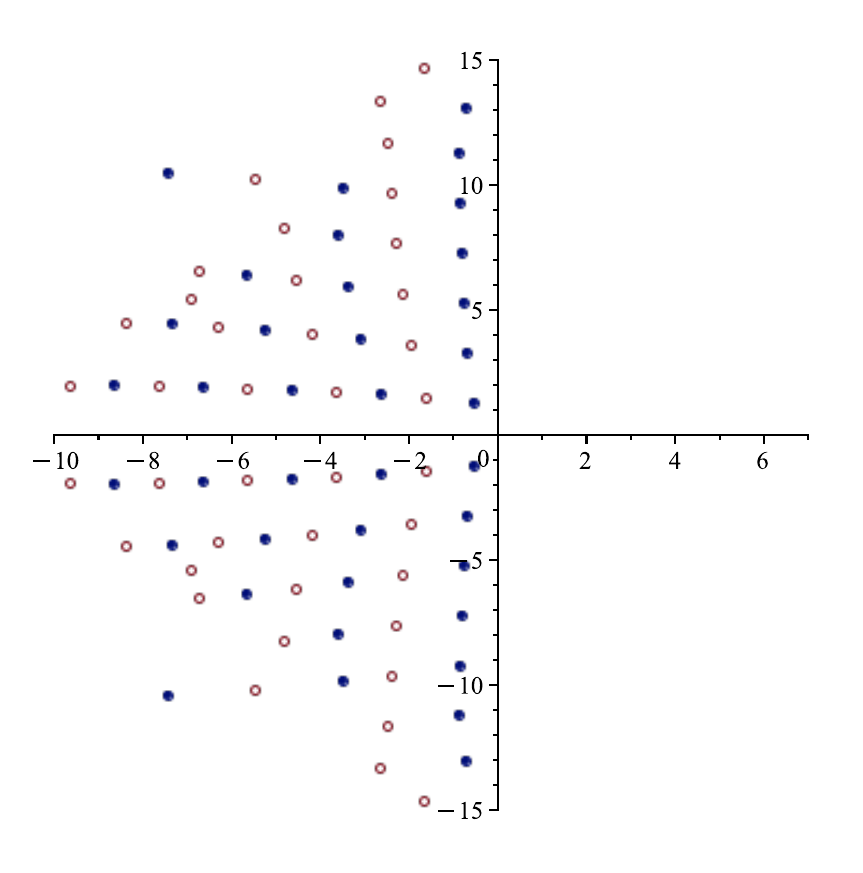}
\caption{The poles of the Pad\'e approximants $[60,40]$ for $P$ (circles) and $Q$ (solid circles).}
\label{fig:pade40}
\end{figure}

Next, we computed the Pad\'e approximants $[60/40]$ for the functions $P$ and $Q$. The poles of these approximants are plotted in Figure \ref{fig:pade40} with
circles for the poles of the Pad\'e approximant for $P$ and solid circles for the poles of the Pad\'e approximant of $Q$. Two large real poles for the latter (one positive, one negative) are not in the picture. The theorem of de Montessus de Ballore only guarantees that the poles closest to the origin are approximating the poles of $P$ and $Q$, and
hence the poles of $r$. Observe that the poles of $P$ and $Q$ have an interlacing pattern.

\subsection{Pad\'e approximation near $\infty$}    \label{sec:PAinf}
We tried to use Pad\'e approximation near infinity with $(n=m)$ using the series in Section \ref{sec:2.4}.
The poles turned out to be real and came in pairs $(-z_k,z_k)$. This probably happens because there are infinitely many poles of $r$ in the complex plane and they are unbounded,
so that it not possible to tell which ones are the closest to infinity, and the Pad\'e approximant generates only real zeros, symmetric around the origin. 
The series expansions in Section \ref{sec:2.4}  hold on the positive real line and hence give little information of what happens in the half plane $\Re z < 0$.
Note that $r$ does not have poles on the positive real axis, so that the poles of these Pad\'e approximants are not of use in approximating the poles of $r$. 

\section{Hermite-Pad\'e approximation}   \label{sec:HermitePade}

\subsection{Hermite-Pad\'e (type I) approximation near $0$}
\begin{figure}[ht]
\centering
\includegraphics[width=4in]{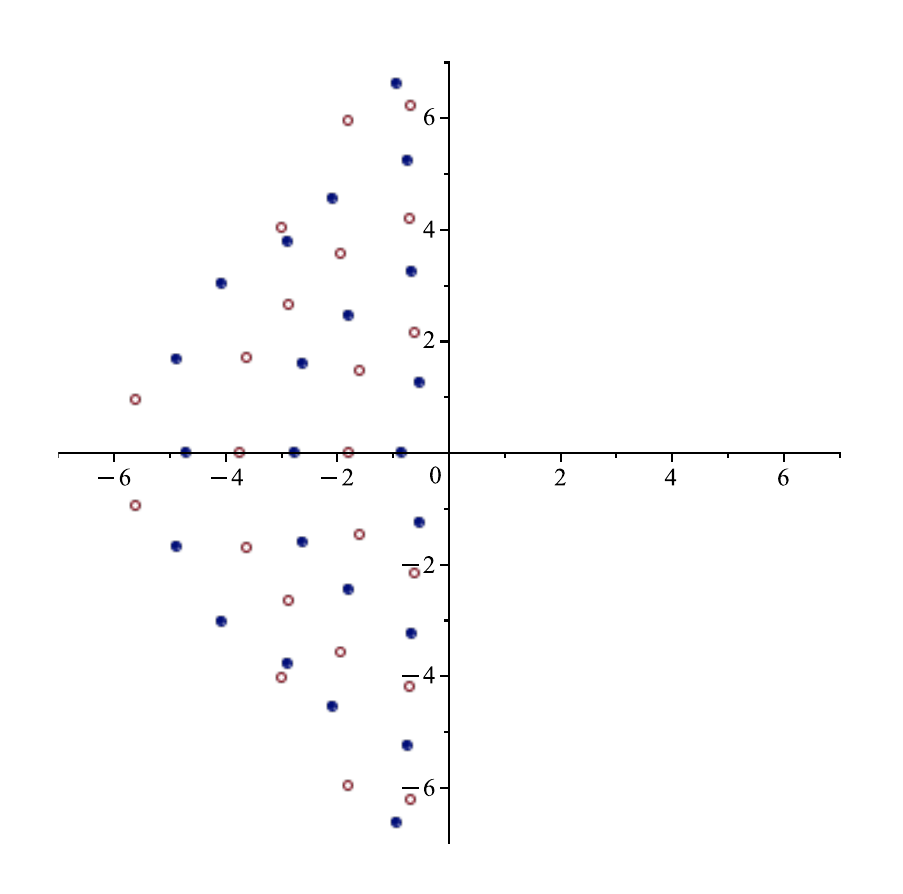}
\caption{The zeros of the type I  Hermite-Pad\'e approximant $[25,25]$ (circles for $A_{25}$ and solid circles for $B_{25}$).}
\label{fig:HP25}
\end{figure}

Next we use Hermite-Pad\'e approximation to the pair $(P,Q)$ near zero, i.e., we look for polynomials $A_n$ and $B_m$, with $\deg A_n \leq n$ and $\deg B_m \leq m$, such that
\[    A_n(x) P(x) + B_m(x) Q(x) = \mathcal{O}(x^{n+m+1}), \qquad x \to 0.  \]
The polynomial $A_n$ tries to cancel the poles of $P$ and $B_m$ tries to cancel the poles of $Q$. The poles of $R=2PQ/\pi$ are then a combination of the poles of $P$ and $Q$.
We have computed the Hermite-Pad\'e approximation with $n=m=25$ and plotted the zeros of $A_{25}$ (which approximate the poles of $p$) as circles in Figure \ref{fig:HP25},
and the zeros of $B_{25}$ as solid circles. The poles in Figure \ref{fig:pade25} seem to be a subset of the points in Figure \ref{fig:HP25}. This means that $A_n$ and $B_n$ contain more information and not all their zeros correspond to poles of $p$ and $q$. Some of the zeros of $A_n$ correspond to zeros of $q$ and some of the zeros of $B_n$ correspond to zeros of $p$.

\subsection{Hermite-Pad\'e (type II) approximation near $0$}
We have also computed simultaneous (type II) Hermite-Pad\'e approximants for the functions $P$ and $Q$. In particular we 
looked for a polynomial $P_{n,m}$ of degree $\leq n+m$, and polynomials $Q_{N-n}$ of degree $\leq N-n$ and $R_{N-m}$ of degree $\leq N-m$ for which 
\begin{eqnarray}
     P_{n,m}(x) P(x) - Q_{N-n}(x) &=& \mathcal{O}(x^{N+1}), \\
     P_{n,m}(x) Q(x) - R_{N-m}(x) &=& \mathcal{O}(x^{N+1}),  
\end{eqnarray}
as $x \to 0$. Here $N$ is  an integer $\geq \max(n,m)$.The rational functions $Q_{N-n}(x)/P_{n,m}(x)$ and $R_{N-m}(x)/P_{n,m}(x)$ are rational approximants for
$P$ and $Q$ respectively, with a common denominator. It is expected that the zeros of this denominator approximate the poles
of $P$ and $Q$, and hence the poles of the function $r$.
There is a  de Montessus de Ballore theorem for type II Hermite-Pad\'e approximants
that allows to conclude that for fixed $n,m$  the zeros of $P_{n,m}$ converge to the poles of $P$ and $Q$ (and thus to the poles of $r$) when $N \to \infty$, i.e.,
the degrees of the numerators $Q_{N-n}$ and $R_{N-m}$ tend to infinity (see Graves-Morris and Saff \cite{GM-Saff}, Cacoq, de la Calle Ysern and L\'opez Lagomasino \cite{CYL1,CYL2}).

\begin{theorem}[Thm. 1.4 in \cite{CYL2}]
Suppose $(P,Q)$ has exactly $n+m$ system poles $z_k$ with respect to $(n,m)$ in $|z| < R$ , then the common denominator $P_{n,m}$ is uniquely determined for all sufficienly large $N$ and there exists a polynomial $\pi_{m+n}$ of degree $m+n$ such that
\[ \limsup_{N \to \infty} \|  P_{n,m} - \pi_{n+m}\|^{1/N}  = \frac{\max_{1 \leq k \leq n+m} |z_k|}{R} < 1.  \] 
\end{theorem}

\begin{figure}[ht]
\centering
\includegraphics[width=4in]{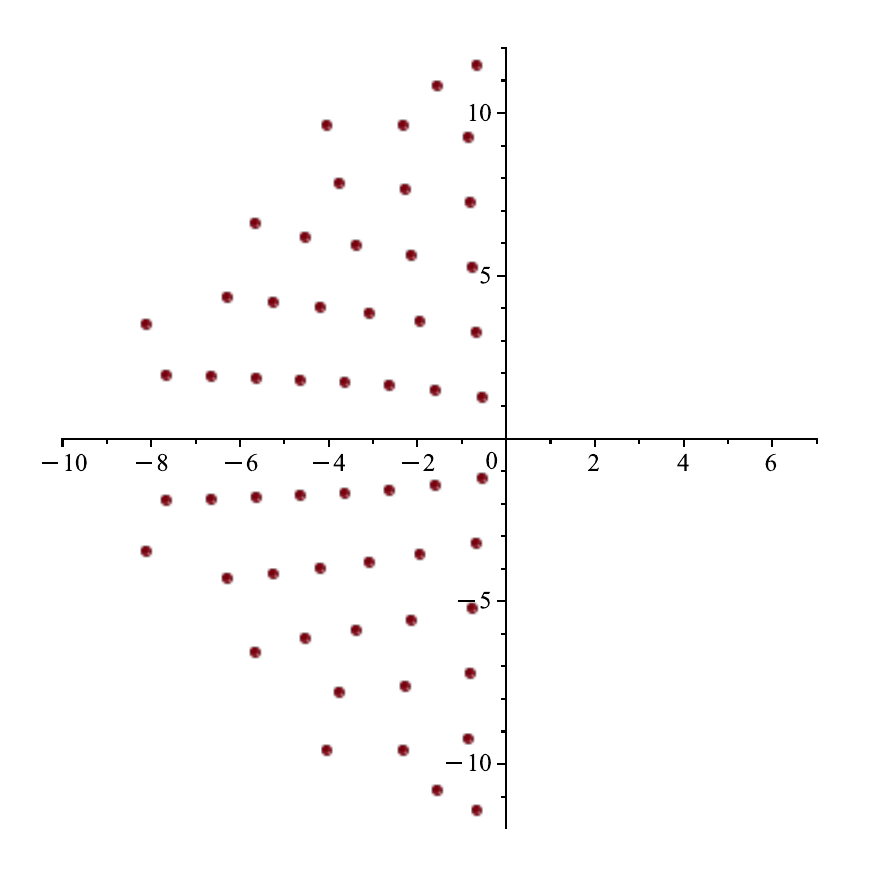}
\caption{Zeros of the common denominator of type II Hermite-Pad\'e approximation with $n=m=30$ and $N=100$.}
\label{fig:HP60}
\end{figure}

We refer to \cite{CYL2} for more details and the notion of system poles.  We believe that the poles of $P$ and $Q$ are system poles for large enough $(n,m)$ because the poles are 
well separated and there are no no other singularities. 
In Figure \ref{fig:HP60} we plotted the zeros of the common denominator $P_{n,m}$ for $n=m=30$ and $N=100$, but only those that are
in the square $[-10,7] \times [-12,12]$. Outside this square there is one large negative real zero (around $-32.48$), one positive real zero (around $9.64$) and a complex pair with
positive real part (approximately $0.42 \pm 13.15 i$). These zeros are not approximating poles of $P$ or $Q$, but of the remaining $56$ zeros there are $54$ that form a nice pattern  and are approximating poles of $r$. The fact that not all $60$ zeros are approximating the poles of $r$ is because $N$ is not sufficiently large, compared to $n+m$.
 
\section{Two-point Pad\'e approximation}   \label{sec:twopoint}
Having a series expansion near $0$ and near $\infty$, it is of interest to use a two-point Pad\'e approximant near $0$ and $\infty$ for the function $r$.
This means that we look for polynomials $P_n$ and $Q_n$ of degree $\leq n$ such that
\begin{eqnarray}
    P_n(x) r(x) - Q_n(x) &=& \mathcal{O}(x^N), \qquad x \to 0,    \label{TPA0} \\
    P_n(x) r(x) - Q_n(x) &=& \mathcal{O}(x^{N-n-1}), \qquad x \to \infty,   \label{TPAinf}
\end{eqnarray}
where $N$ is a positive integer with $N \leq 2n$ that determines how many interpolation conditions are used at $0$ and at $\infty$: the condition
\eqref{TPA0} gives $N$ conditions and condition \eqref{TPAinf} gives $2n-N+1$ conditions. This means that we have $2n+1$ linear and homogeneous equations
for the $2n+2$ unknown coefficients of $P_n$ and $Q_n$. 
See, for example \cite{McCabe-Murphy}.
If we change $x$ to $1/x$ then the condition \eqref{TPAinf} becomes
\[     P_n(1/x) r(1/x) - Q_n(1/x) = \mathcal{O}(x^{n-N+1}), \qquad x \to 0. \]
If we denote $P_n^*(x) = x^n P_n(1/x)$ and $Q_n^*(x)=x^n Q_n(1/x)$, then the coefficients of $P_n^*$ and $Q_n^*$ are those of $P_n$ and $Q_n$ but in reversed order, and
the interpolation condition \eqref{TPAinf} becomes
\begin{equation}   \label{TPA*}
     P_n^*(x) r(1/x) - Q_n^*(x) = \mathcal{O}(x^{2n-N+1}), \qquad x \to 0.
\end{equation}      
The two-point Pad\'e approximant thus corresponds to solving two Pad\'e-type approximation problems at $0$. 

 \begin{figure}[ht]
\centering
\includegraphics[width=3in]{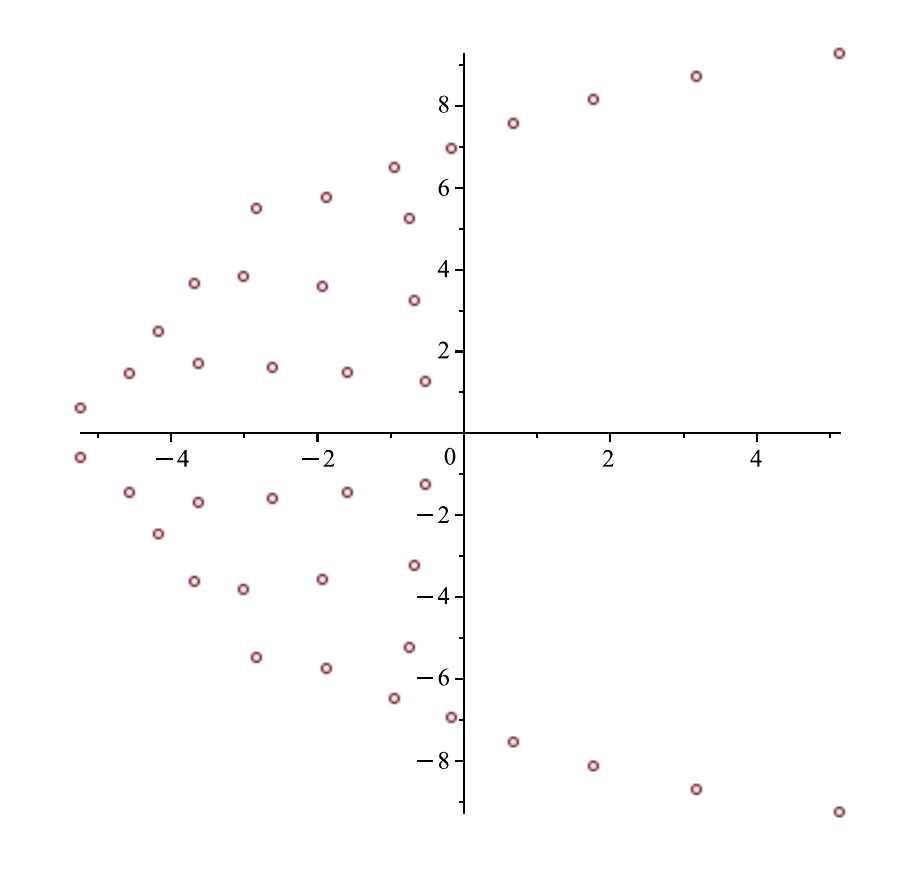}
\includegraphics[width=3in]{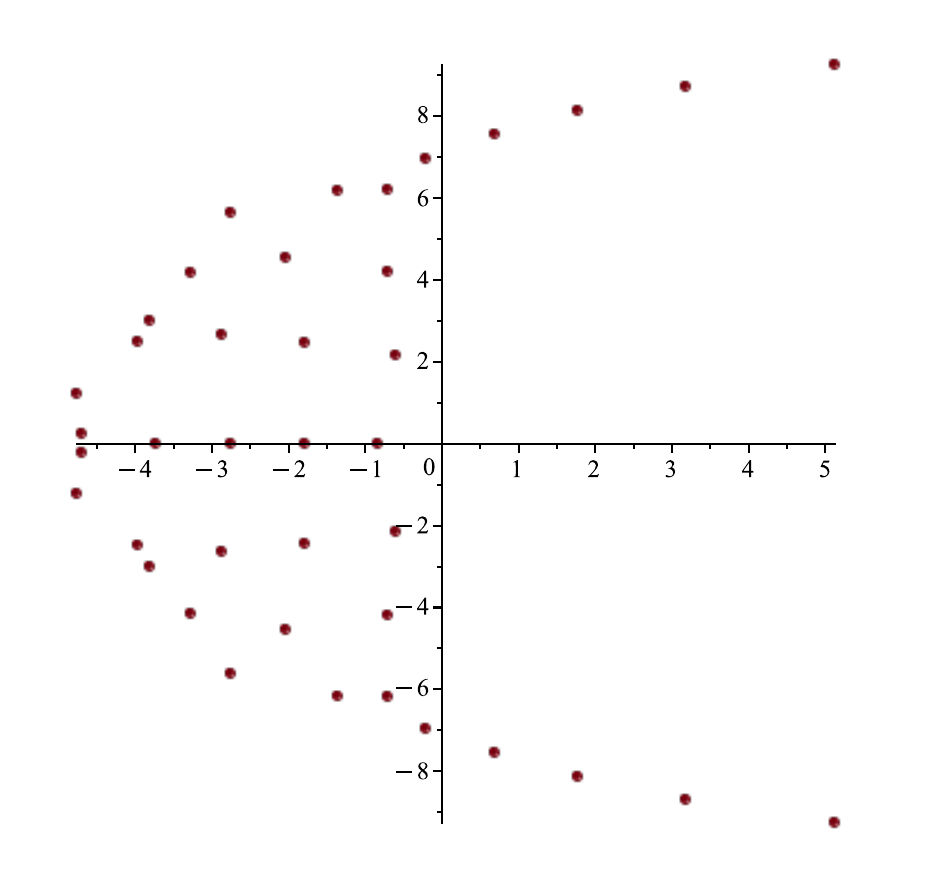}
\caption{Zeros of $P_n$ (left) and $Q_n$ (right) of the two-point Pad\'e approximant with $n=40$ and $N=60$.}
\label{fig:TPA}
\end{figure}

In Figure \ref{fig:TPA} we plotted the zeros of $P_n$ (left) and $Q_n$ (right) for $n=40$ and $N=60$.
Note that $P_n$ has poles with positive real part, which we don't expect, but $Q_n$ also has zeros at those points, so that these unexpected poles are cancelled in the ratio $Q_n/P_n$ which approximates the function $r$. There are fewer unexpected poles if we increase the number $N$ of interpolation conditions at $0$. For $N=2n$ we are back at
the Pad\'e approximant at $0$ from Section \ref{sec:PA0} with $m=n$, and for $N=0$ we get the Pad\'e approximant at infinity from Section \ref{sec:PAinf}. 
 
 \section{Conclusion}
 We have used Pad\'e approximation to produce approximations of the poles of the Painlev\'e transcendent $r$ appearing in the expression of the
 Fredholm determinant of the sine kernel (Theorem \ref{thm:Fsin}).
 The Pad\'e approximant  $[m/n]$ near $0$ for the function $r$ has the property that if we fix the degree of the denominator $n$ and let $m \to \infty$, then the zeros
 of the denominator converge to the first $n$ poles of $r$, i.e, the $n$ poles closest to the origin. This is a consequence of the theorem of de Montessus de Ballore.
We found it more efficient to compute the $[m/n]$ Pad\'e approximants of two functions $P$ and $Q$ near $0$ and then to combine the zeros of the denominators
to find approximations for the first $2n$ poles of $r$. This follows because $r(x)=2P(x)Q(x)/\pi$, where $P(x)=p(x)$ and $Q(x)=q(x)/x$. We observed that 
the poles of $P$ and $Q$ all are in the half plane $\Re z < 0$ and that they seem to have an interlacing behavior (see Figure \ref{fig:pade40}). 

We tried Pad\'e approximation near infinity, but the theorem of de Montessus de Ballore can't be used near infinity since $r$ has infinitely many poles which are accumulating at infinity. Furthermore the series espansions in Section \ref{sec:2.4} are valid on the positive real axis and are no longer true in the half plane $\Re z <0$ where the poles are.
The zeros of the denominator of the Pad\'e approximant near infinity are real and are not giving any information about the poles of $r$ in the complex plane.

We also tried using Hermite-Pad\'e approximation (simultaneous Pad\'e approximation) to the pair $(P,Q)$ near $0$. The type I Hermite-Pad\'e approximant near $0$
tries to make $A_n(x)P(x)+B_n(x)Q(x)$ as small as possible near $0$, This means that near poles of $P$ one needs zeros of $A_n$ and near poles of $Q$ one needs
zeros of $B_n$ to keeps the linear form bounded. But at zeros of $P$ one needs zeros of $B_n$ and at zeros of $Q$ one needs zeros of $A_n$ to make sure that
the linear form has a large zero at the origin. This means that some of the zeros of $A_n$ approximate poles of $P$ and some approximate zeros of $Q$, while some
of the zeros of $B_n$ approximate poles of $Q$ and some approximate zeros of $P$. This makes it difficult to determine which zeros of $A_n$ and $B_n$ are approximating the poles of $r$. Type II Hermite-Pad\'e approximation turned out to be more useful. One can use the zeros of the common denominator $P_{n,m}$ as
approximants for the poles of $P$ and $Q$ because these poles are system poles of $(P,Q)$, see \cite{BosLopez}. 
The zeros of the common denominator $P_{n,m}$  keep the product $P_{n,m}P$ and
$P_{n,m}Q$ bounded near the poles of $P$ and $Q$ and the zeros of $P_{n,m}$ are therefore approximations for $n+m$ poles of $r$ nearest to the origin (see Figure \ref{fig:HP60}). There is an analogue of the theorem of de Montessus de Ballore that allows to conclude that the zeros of $P_{n,m}$ converge to the poles of $P$ and
$Q$ if the degree of the numerators tends to infinity. 
 
 Finally we also tried two-point Pad\'e approximation near $0$ and infinity. We found that the information near infinity was not useful and in fact resulted in
 getting zeros of the denominator at unexpected places, namely with positive real part. These poles do not approximate any poles of $r$, and the two-point
 Pad\'e approximant compensates this by having zeros  near those poles, so that the zeros of the two-point Pad\'e approximant cancel the poles of the two-point Pad\'e approximant. The two-point Pad\'e approximant therefore was not so useful for approximating the poles of $r$, but the rational approximant $Q_n/P_n$ gives a good approximation of the function $r$ on $[0,\infty)$ with the correct behavior near $0$ and $\infty$.
 
In the present work the approximations for the poles use Pad\'e approximation of high order which requires to solve a large linear system of equations. Furthermore the poles
are always zeros of a polynomial of large degree, and hence one needs to have a way to compute zeros of a large degree polynomial.
We have done this using computeralgebra (Maple 2026) with high precision (\texttt{Digits:=100}). Furthermore we can only approximate poles close to the origin
since we are using Pad\'e approximants near $0$. Fornberg and Weideman \cite{FW} have developed a numerical method solving the Painlev\'e differential equation
using a low order Pad\'e approximant around $x_k$ to approximate $r(x_{k+1})$ for computing $r$ at points $x_k$ on a path in the complex plane.
It would be nice to see how their method approximates poles for the Painlev\'e transcendent $r$.

 \section*{Acknowledgments}
 No AI was used to produce this paper.
 Maple 2026 was used for the computations in high precision (\texttt{Digits:=100}) and to produce the figures. 
 I am very grateful to George Labahn and Berhard Beckermann for useful comments and conversations.

\end{document}